\documentclass[submission]{FPSAC2021}

\usepackage{tikz}
\usetikzlibrary{positioning}

\theoremstyle{plain}
\newtheorem{Th}{Theorem}[section]

\newtheorem{Prop}[Th]{Proposition}

\newcommand{\St}{\operatorname{St}}
\newcommand{\type}{\operatorname{type}}

\newcommand{\Sym}{\mathfrak{S}}

\newcommand{\length}{\operatorname{length}}

\newcommand{\code}{\operatorname{code}}
\newcommand{\Z}{\mathbb{Z}}

\newcommand{\x}{\mathbf{x}}

 \theoremstyle{definition}
\newtheorem{Def}[Th]{Definition}

\newtheorem{Rem}[Th]{Remark}
\newtheorem{?}[Th]{Problem}
\newtheorem{Ex}[Th]{Example}

\title[Schubert polynomials and TASEP]{Schubert polynomials and the inhomogeneous TASEP on a ring}

\author[Donghyun Kim and Lauren Williams]{Donghyun Kim\thanks{\href{mailto:donghyun_kim@berkeley.edu}{donghyun\_kim@berkeley.edu}.}\addressmark{1}, \and Lauren K. Williams\thanks{L.W. was partially supported by NSF grant DMS-1854512}\addressmark{2}}

\address{\addressmark{1}Department of Mathematics, University of California at Berkeley, CA \\ \addressmark{2}Department of Mathematics, Harvard University, Cambridge, MA}

\received{\today}

\abstract{Consider a lattice of n sites arranged
around a ring, with the $n$ sites occupied by 
particles of weights $\{1,2,\dots,n\}$; the possible
arrangements of particles in sites thus corresponds to the 
$n!$ permutations in $S_n$.
The \emph{inhomogeneous totally asymmetric simple
exclusion process} (or TASEP) is a Markov chain 
on the set of permutations, 
in which two adjacent particles of weights $i<j$
swap places at rate $x_i - y_{n+1-j}$ if the particle of weight $j$ is to the right of the particle of weight $i$.  (Otherwise nothing happens.)  
In the case that $y_i=0$ for all $i$, the stationary
distribution was conjecturally linked to Schubert polynomials
by Lam-Williams, and explicit formulas for steady 
state probabilities were subsequently given in terms
of multiline queues by Ayyer-Linusson and Arita-Mallick.
In the case of general $y_i$,  Cantini showed that $n$ of the 
$n!$ states have probabilities proportional to double
Schubert polynomials.  In this paper
we introduce the class of \emph{evil-avoiding permutations},
which are the permutations avoiding the patterns
$2413, 4132, 4213$ and $3214$.
We show that 
there are $\frac{(2+\sqrt{2})^{n-1}+(2-\sqrt{2})^{n-1}}{2}$
evil-avoiding permutations in $S_n$, and for each 
evil-avoiding permutation $w$, we give an explicit formula
for the steady 
state probability $\psi_w$ as a product
of double Schubert polynomials. We also show that the Schubert polynomials that arise in these formulas are flagged Schur functions, and give a bijection in this case between 
multiline queues and semistandard Young tableaux.  }

\keywords{Schubert polynomials, TASEP, multiline queues}

\usepackage[backend=bibtex]{biblatex}
\begin{document}

\maketitle

\section{Introduction}
In recent years, there has been a lot of work 
on interacting particle models such as the 
\emph{asymmetric simple exclusion process} (ASEP),
a model in which particles hop on a one-dimensional
lattice subject to the condition that at most one particle
may occupy a given site.  The ASEP on a 
one-dimensional lattice with open boundaries has 
been linked to Askey-Wilson polynomials and Koornwinder
polynomials \cite{ CW1, C2, CW2}, 
while the ASEP on a ring has been linked to 
Macdonald polynomials \cite{CGW, CMW}.
The \emph{inhomogeneous totally asymmetric simple exclusion 
process} (TASEP) is a variant of the exclusion process 
on the ring in which the hopping rate depends on the weight of 
the particles.  In this paper we build on works of 
Lam-Williams \cite{LW}, Ayyer-Linusson \cite{AL},
and especially Cantini \cite{C} to give formulas for 
many steady state probabilities of the inhomogeneous TASEP on a ring in terms of Schubert polynomials.

\begin{Def}\label{def:TASEP}
Consider a lattice with $n$ sites arranged in a ring.
Let $\St(n)$ denote the $n!$ labelings
of the lattice by distinct numbers $1,2,\dots,n$, where each number $i$ is called a 
\emph{particle of weight $i$}.  
The \emph{inhomogeneous TASEP on a ring of size $n$} is a Markov chain with state space $\St(n)$
where at each time $t$ a swap of two adjacent particles may occur: a particle of weight $i$ on the left swaps its position with a particle of weight $j$ on the right with transition rate $r_{i,j}$ given by:
$$
r_{i,j} = 
\begin{cases}  
x_i-y_{n+1-j} \text{ if $i<j$} \\
0 \text{ otherwise.}
\end{cases}
$$
\end{Def}

\begin{figure}[h]
\begin{center}

\begin{tikzpicture}[scale=1.5]
\node (0) [label={[red]below: $x_1+x_2-y_1-y_2$}] at (0,0) {$321$};
\node (12)  [label={[red]right: $x_1+x_2-y_1-y_2$}] at (1,2) {$213$};
\node (21)  [label={[red]left: $x_1+x_2-y_1-y_2$}] at (-1,2) {$132$};
\node (1) [label={[red]left: $x_1-y_1$}] at (-1,1) {$312$};
\node (2) [label={[red]right: $x_1-y_1$}] at (1,1) {$231$};
\node (121) [label={[red]above: $x_1-y_1$}] at (0,3) {$123$};

\draw[->] (0) -- node[near start,blue]{$x_1-y_1$} (121);
\draw[->](1)--node[near start,sloped, above=0.5pt,blue]{}(12);
\draw[->](2)--node[near start,sloped, above=0.5pt,blue]{} (21);
\draw[->](121)--node[left=0.5pt,blue]{$x_2-y_1$}(21);
\draw[->](121)--node[right=0.5pt,blue]{$x_1-y_2$}(12);
\draw[->](21)--node[left=0.5pt,blue]{$x_1-y_1$}(1);
\draw[->](12)--node[right=0.5pt,blue]{$x_1-y_1$}(2);
\draw[->](1) -- node[left=0.5pt,blue]{$x_1-y_2$}(0);
\draw[->](2)--node[right=0.5pt,blue]{$x_2-y_1$}(0);
\end{tikzpicture}
\caption{The state diagram for the inhomogeneous TASEP on $\St(3)$, with transition rates shown in blue, and steady state probabilities $\psi_w$ in red.
Though not shown, the transition rate $312\to 213$ is $x_2-y_1$ and the 
transtition rate $231\to 132$ is $x_1-y_2$.}
\label{fig:S3}
\end{center}
\end{figure}
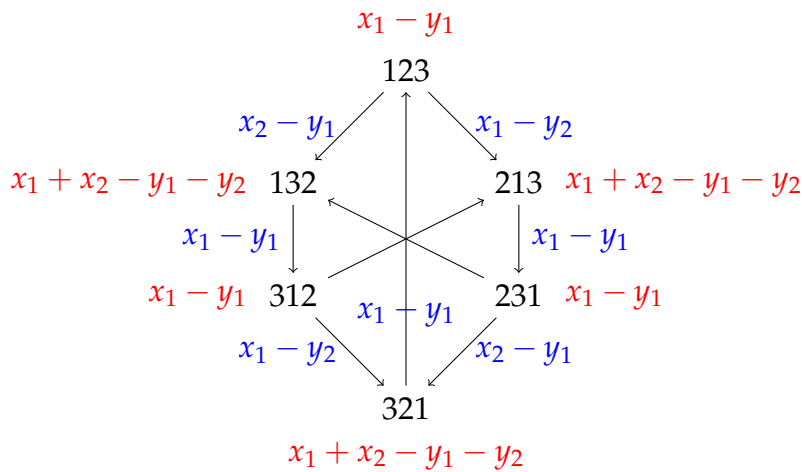

 In what follows, we will identify each state with a permutation in $S_n$.   Following \cite{LW, C}, we multiply all steady state probabilities for $\St(n)$ by the same constant, obtaining 
 ``renormalized'' steady state probabilities  $\psi_w$, so that 
  \begin{equation}\label{eq:normalization}
 \psi_{123 \dots n} = \prod_{i<j} (x_i-y_{n+1-j})^{j-i-1}.
 \end{equation}
 See Figure \ref{fig:S3} for the state diagram 
 when $n=3$.  

In the case that $y_i=0$, Lam and Williams \cite{LW} 
studied this model\footnote{However the convention of \cite{LW} was slightly different; it corresponds to labeling states by the inverse of the permutations we use here.}
and conjectured that after a suitable normalization, each 
steady state probability $\psi_w$ can be written as a monomial
factor times a positive sum of Schubert polynomials, see 
\cref{table:1} and \cref{table:2}.  
They also gave an explicit formula for the monomial
factor, and conjectured that under certain conditions on 
$w$,  $\psi_w$ is a multiple of 
a particular Schubert polynomial. Subsequently  Ayyer and Linusson \cite{AL} gave a conjectural combinatorial
formula for the stationary distribution in terms of \emph{multiline queues}, which was proved by 
Arita and Mallick \cite{AM}.
 In \cite{C}, Cantini introduced the version of the model given in \cref{def:TASEP}\footnote{We
 note that in \cite{C}, the rate $r_{i,j}$  was $x_i-y_{j}$ rather than $x_i-y_{n+1-j}$ as we use in \cref{def:TASEP}.}
 with $y_i$ general, and gave a series of 
 \emph{exchange equations} relating the components of 
 the stationary distribution.  This allowed him to give explicit formulas for the steady state probabilities for $n$ of the 
$n!$ states as products of double
Schubert polynomials. 

\begin{table}[h!]
\centering
\begin{tabular}{|c c |}
    \hline
    State $w$ & Probability $\psi_w$\\
    \hline 
    1234 & $(x_1-y_1)^2 (x_1-y_2)(x_2-y_1)$\\
    1324 & $(x_1-y_1) \Sym_{1432}$\\
    1342 & $(x_1-y_1)(x_2-y_1) \Sym_{1423}$\\
    1423 & $(x_1-y_1)(x_1-y_2)(x_2-y_1) \Sym_{1243}$\\
    1243 & $(x_1-y_2)(x_1-y_1) \Sym_{1342}$\\
    1432 & $\Sym_{1423} \Sym_{1342}$\\ 
    \hline
    \end{tabular}
    \caption{The renormalized steady state probabilities for $n=4$.}
    \label{table:1}
    \end{table}

In this paper we build on \cite{C, AL, AM}, and give
many more explicit formulas for steady state probabilities
in terms of Schubert polynomials: in particular, we give
a formula for $\psi_w$ as a product of (double) Schubert polynomials whenever $w$ is \emph{evil-avoiding}, that is, it
avoids the patterns
$2413, 4132, 4213$ and $3214$.\footnote{We call these permutations \emph{evil-avoiding} because if one 
replaces $i$ by $1$, $e$ by $2$, $l$ by $3$, and $v$ by $4$, 
then \emph{evil} and its anagrams \emph{vile, veil} and \emph{leiv} become the four patterns $2413, 4132, 4213$ and $3214$.  Note that Leiv is a name of Norwegian origin meaning ``heir.''}  
We show that 
there are $\frac{(2+\sqrt{2})^{n-1}+(2-\sqrt{2})^{n-1}}{2}$
evil-avoiding permutations in $S_n$, so this gives a 
substantial generalization of Cantini's previous result
\cite{C} 
in this direction.  We also prove the 
monomial factor conjecture from \cite{LW}.  Finally, we show that the Schubert polynomials that arise in our formulas are flagged Schur functions, and give a bijection in this case between multiline queues and semistandard Young tableaux.

\begin{center}
\begin{table}[h!]
\centering
    \begin{tabular}{|c c |}
    \hline
    State $w$ & Probability $\psi_w$\\
    \hline 
    12345 & $\mathbf{x}^{(6,3,1)}$\\
    12354 & $\mathbf{x}^{(5,2,0)} \Sym_{13452}$\\
    12435 & $\mathbf{x}^{(4,1,0)} \Sym_{14532}$\\
    12453 & $\mathbf{x}^{(4,1,1)} \Sym_{14523}$\\
    12534 & $\mathbf{x}^{(5,2,1)} \Sym_{12453}$\\
    12543 & $\mathbf{x}^{(3,0,0)} \Sym_{14523} \Sym_{13452}$\\
    13245 & $\mathbf{x}^{(3,1,1)} \Sym_{15423}$\\
    13254 & $\mathbf{x}^{(2,0,0)} \Sym_{15423} \Sym_{13452}$\\
    13425 & $\mathbf{x}^{(3,2,1)} \Sym_{15243}$\\
    13452 & $\mathbf{x}^{(3,3,1)} \Sym_{15234}$\\
    13524 & $\mathbf{x}^{(2,1,0)} (\Sym_{164325}+\Sym_{25431})$\\
    13542 & $\mathbf{x}^{(2,2,0)} \Sym_{15234} \Sym_{13452}$\\
    14235 & $\mathbf{x}^{(4,2,0)} \Sym_{13542}$\\
    14253 & $\mathbf{x}^{(4,2,1)} \Sym_{12543}$\\
    14325 & $\mathbf{x}^{(1,0,0)} 
(\Sym_{1753246}+\Sym_{265314}+\Sym_{2743156}+\Sym_{356214}+\Sym_{364215}+\Sym_{365124})$\\
    14352 & $\mathbf{x}^{(1,1,0)} \Sym_{15234} \Sym_{14532}$\\
    14523 & $\mathbf{x}^{(4,3,1)} \Sym_{12534}$\\
    14532 & $\mathbf{x}^{(1,1,1)} \Sym_{15234} \Sym_{14523}$\\
    15234 & $\mathbf{x}^{(5,3,1)} \Sym_{12354}$\\
    15243 & $\mathbf{x}^{(3,1,0)}(\Sym_{146325}+\Sym_{24531})$\\
    15324 & $\mathbf{x}^{(2,1,1)} (\Sym_{15432}+\Sym_{164235})$\\
    15342 & $\mathbf{x}^{(2,2,1)} \Sym_{15234}\Sym_{12453}$\\
    15423 & $\mathbf{x}^{(3,2,0)} \Sym_{12534}\Sym_{13452}$\\
    15432 & $\Sym_{15234} \Sym_{14523}\Sym_{13452}$\\
    \hline\end{tabular}
    \caption{The renormalized steady state probabilities for $n=5$, when
    each $y_i=0$.  In the table, $\mathbf{x}^{(a,b,c)}$ denotes $x_1^a x_2^b x_3^c$.}
    \label{table:2}
    \end{table}
\end{center}

In order to state our main results, we need a few definitions.
First, we say that two states $w$ and $w'$ are \emph{equivalent}, and write
$w\sim w'$, if one state is a cyclic shift of the other, e.g.   
$(w_1,\dots,w_n) \sim (w_2,\dots,w_n,w_1)$.  Because of the cyclic symmetry inherent in the definition of the TASEP on a ring, it is clear that the probabilities of states $w$ and $w'$ are equal whenever $w\sim w'$.  
We will therefore often assume, without loss of generality, 
that $w_1=1$.  Note that up to cyclic shift, $\St(n)$ contains $(n-1)!$ states.

\begin{Def}\label{def:kGrassmannian}
Let $w=(w_1,\dots,w_n)\in \St(n)$. 
We say that $w$ is a 
\emph{$k$-Grassmannian permutation}, 
and we write
$w\in \St(n,k)$ if: $w_1=1$; $w$ is \emph{evil-avoiding}, i.e. 
    $w$ avoids
    the patterns $2413$, $3214$, $4132$, and $4213$; and
 $w^{-1}$ has exactly $k$ \emph{descents}, equivalently, 
    there are exactly $k$ letters $a$ in $w$ such that $a+1$ appears
    to the left of $a$ in $w$.
\end{Def}

\begin{Def}\label{def:bijpartition}
We associate to each $w\in\St(n,k)$ a sequence of partitions
$\Psi(w)=(\lambda^1,\dots,\lambda^k)$ 
 as follows. Write the Lehmer code of $w^{-1}$ as 
 $\code(w^{-1}) = c = (c_1,\dots,c_n);$
since $w^{-1}$ has $k$ descents, 
$c$ has $k$ \emph{descents} in positions we denote by
$a_1,\dots,a_k$.  We also set $a_0=0$.  For $1 \leq i \leq k$,
we define
$\lambda^i = (n-{a_i})^{a_i} - (\underbrace{0,\cdots,0}_\text{$a_{i-1}$},c_{a_{i-1}+1},
c_{a_{i-1}+2},\dots, c_{a_i}).$
\end{Def}
See \cref{table:3} for examples of the map $\Psi(w)$.

\begin{Def}\label{def:g} Given a positive integer $n$ and a partition $\lambda$ of length $\leq(n-2)$, we define an integer vector 
$g_n(\lambda)=(v_1,\dots,v_n)$ of length $n$ as follows.  Write $\lambda=(\mu_1^{k_1},\cdots,\mu_l^{k_l})$ where $k_i>0$ and $\mu_1>\cdots>\mu_l$.
We assign values to the entries $(v_1,\dots,v_n)$ by performing the following step for $i$ from $1$ to $l$.
\begin{itemize}
    \item (Step $i$) Set $v_{n-\mu_i}$ equal to $\mu_i$. Moving to the left, assign the value $\mu_i$ to the first $(k_i-1)$ unassigned components.
\end{itemize}
After performing Step $l$, we assign the value $0$ to any entry $v_j$ which has not yet been given a value.
\end{Def}

Note that in Step 1, we set $v_{n-\mu_1},v_{n-\mu_1-1},\cdots,v_{n-\mu_1-k_1+1}$ equal to $\mu_1$.

\begin{Ex}
\begin{align*}
    g_5((2,1,1))=(0,1,2,1,0)\\
    g_6((3,2,2,1))=(0,2,3,2,1,0)\\
     g_6((3,1,1))=(0,0,3,1,1,0).
\end{align*}
\end{Ex}

The main result of this paper is
\cref{thm:main1}. We state here our main result in the case that each $y_i=0$.  The definition of Schubert polynomial 
can be found in \cref{sec:background}.

\begin{Th}\label{thm:main0}
Let $w\in \St(n,k)$ be a 
\emph{$k$-Grassmannian} permutation, as 
in \cref{def:kGrassmannian}, and   
let $\Psi(w)=(\lambda^1,\dots,\lambda^k)$.
 Adding trailing $0$'s if necessary, we view each 
partition $\lambda^i$ as 
a vector in $\Z_{\geq 0}^{n-2}$, and set 
$\mu:=(\binom{n-1}{2}, \binom{n-2}{2},\dots ,\binom{2}{2})-\sum_{i=1}^k \lambda^i.$  Then when each $y_i=0$,
the renormalized steady state probability 
$\psi_w$ is given by 
$$\psi_w = \mathbf{x}^{\mu} \prod_{i=1}^k \Sym_{g_n(\lambda^i)},$$
where $\Sym_{g_n(\lambda^i)}$ is the Schubert polynomial
associated to the permutation with Lehmer code
$g_n(\lambda^i)$, and $g_n$ is given by \cref{def:g}.

Equivalently, writing $\lambda^i = (\lambda_1^i, \lambda_2^i,\dots)$, we have that 
$$\psi_w = \mathbf{x}^{\mu} \prod_{i=1}^k
s_{\lambda^i}(X_{n-\lambda_1^i}, X_{n-\lambda_2^i},\dots),$$
where $s_{\lambda^i}(X_{n-\lambda_1^i}, X_{n-\lambda_2^i},\dots)$ denotes the flagged
Schur polynomial associated to shape $\lambda^i$,
where the semistandard tableaux entries in row $j$ are bounded above by 
$n-\lambda_j^i$.  
\end{Th}

We illustrate \cref{thm:main0} in \cref{table:3} in the case that $n=5$.

\begin{center}
\begin{table}[h!]
\centering
    \begin{tabular}{|c c c c |}
    \hline
k &     $w\in \St(5,k)$ & $\Psi(w)$ & probability $\psi_w$ \\
    \hline 
0 &    12345 & $\emptyset$ & $\x^{(6,3,1)}$\\
\hline 
 1&   12354 & $(1,1,1)$ & $\x^{(5,2,0)} \Sym_{13452}$ \\
   1& 12435 & $(2,2,1)$ & $\x^{(4,1,0)} \Sym_{14532}$\\
    1& 12453 & $(2,2)$ & $\x^{(4,1,1)} \Sym_{14523}$\\
    1& 12534 & $(1,1)$ & $\x^{(5,2,1)} \Sym_{12453}$\\
    1& 13245 & $(3,2)$ & $\x^{(3,1,1)} \Sym_{15423}$\\
    1& 13425 & $(3,1)$& $ \x^{(3,2,1)} \Sym_{15243}$\\
    1& 13452 & $(3)$ & $\x^{ (3,3,1)} \Sym_{15234}$ \\
    1& 14235 & $(2,1,1)$ & $\x^{(4,2,0)}  \Sym_{13542}$\\
    1& 14253 & $(2,1)$ & $ \x^{(4,2,1)} \Sym_{12543}$\\
    1& 14523 & $(2)$ & $ \x^{(4,3,1)} \Sym_{12534}$ \\
    1& 15234 & $(1)$ & $\x^{(5,3,1)} \Sym_{12354}$ \\
        \hline 
    2& 12543 & $(2,2), (1,1,1)$& $\x^{(3,0,0)} \Sym_{14523} \Sym_{13452}$ \\
    2& 13254 & $(3,2), (1,1,1)$ & $\x^{(2,0,0)} \Sym_{15423} \Sym_{13452}$ \\
    2& 13542 & $(3), (1,1,1)$ & $\x^{(2,2,0)} \Sym_{15234} \Sym_{13452}$\\
    2& 14352 & $(3), (2,2,1)$ & $\x^{(1,1,0)} \Sym_{15234} \Sym_{14532}$ \\  
    2& 14532 & $(3), (2,2)$ & $\x^{(1,1,1)} \Sym_{15234} \Sym_{14523}$\\
    2& 15342 & $(3), (1,1)$ & $\x^{(2,2,1)} \Sym_{15234}\Sym_{12453}$ \\
    2& 15423 & $(2),(1,1,1)$ & $\x^{(3,2,0)} \Sym_{12534}\Sym_{13452}$ \\
    \hline 
3 &    15432 & $(3), (2,2), (1,1,1)$ & $\Sym_{15234} \Sym_{14523} \Sym_{13452}$ \\
    \hline\end{tabular}
    \caption{Special states $w\in\St(5,k)$ and the corresponding sequences of partitions $\Psi(w)$, together with steady state probabilities $\psi_w$.}
    \label{table:3}
    \end{table}
\end{center}

\begin{Prop}\label{prop:enumeration}
The number of evil-avoiding permutation in 
$S_n$ satisfies the recurrence	
$e(1)=1, e(2)=2, 
e(n)=4 e(n-1)-2e(n-2)$ for $n \geq 3$, and is given 
explicitly as 
\begin{equation}
    \label{eq:a}
e(n)=\frac{(2+\sqrt{2})^{n-1}+(2-\sqrt{2})^{n-1}}{2}.
\end{equation}
This sequence begins as $1, 2, 6, 20, 68, 232$, and 
occurs in
Sloane's encyclopedia as sequence A006012.
The cardinalities 
$|\St(n,k)|$ also occur as sequence 
A331969.
\end{Prop}

\begin{Rem}
Let $w(n,h): = (h, h-1,\dots,2,1,h+1,h+2,\dots,n)\in \St(n).$
In \cite[Corollary 16]{C}, Cantini gives a formula for the steady 
state probability of state $w(n,h)$, as a trivial factor times a product of certain (double) Schubert polynomials.  Note that our main result is a significant generalization of \cite[Corollary 16]{C}.  For example,
for $n=4$, Cantini's result gives a formula for the probabilities of 
three states -- 
$(1,2,3,4), (1,3,4,2),$ and $(1,4,3,2)$.  And for $n=5$, his result gives
a formula for four states -- $(1,2,3,4,5), (1,3,4,5,2),\\ (1,4,5,3,2)$, and $(1,5,4,3,2).$
On the other hand, \cref{thm:main0} gives a formula for all six states when $n=4$ (see \cref{table:1}) and $20$ of the $24$ states when $n=5$.
Asymptotically, since the number of special states in $S_n$ is given by \eqref{eq:a}, \cref{thm:main0} 
gives a formula for roughly $\frac{(2+\sqrt{2})^{n-1}}{2}$ 
out of the 
$(n-1)!$ states of $\St(n).$

Another point worth mentioning is that the Schubert polynomials that occur in the formulas of \cite{C} are all of the form $\Sym_{\sigma(a,n)}$, where 
$\sigma(a,n)$ denotes the permutation $(1, a+1, a+2,\dots, n, 2, 3, \dots,n).$  However, many of the Schubert polynomials arising as (factors) of steady probabilities are not of this form.  Already we see for $n=4$ the Schubert polynomials $\Sym_{1432}$ and $\Sym_{1243}$, which are not of this form.
\end{Rem}

Note that it is common to consider a version of the inhomogeneous TASEP 
in which one allows multiple particles of each 
weight $i$.  This is the version studied in several of the previous references, and also 
in \cite{Mandelshtam} (which primarily considers particles of types $0$, $1$ and $2$).  We plan 
to work in this generality in our subsequent work.  However, since
our focus here is on Schubert polynomials, we restrict to the case of permutations.

\section{Background on permutations and Schubert polynomials}\label{sec:background}

We let $S_n$ denote the symmetric group on $n$ letters, which is a Coxeter group 
generated by the simple reflections $s_1,\dots, s_{n-1}$, where $s_i$ is the simple 
transposition exchanging $i$ and $i+1$.  We let 
$w_0=(n,n-1,\dots, 2,1)$ denote the longest permutation.

For $1 \leq i < n$, we have the \emph{divided difference operator} $\partial_i$ which acts on polynomials 
$P(x_1,\dots,x_n)$ as follows:
$$(\partial_i P)(x_1,\dots,x_n) = \frac{P(\dots, x_i, x_{i+1},\dots ) - P(\dots,x_{i+1},x_i,\dots)}{x_i-x_{i+1}}.$$
If $s_{i_1}\dots s_{i_m}$ is a reduced expression for a permutation $w$, then
$\partial_{i_1} \dots \partial_{i_m}$ depends only on $w$, so we denote this operator
by $\partial_w$.  

\begin{Def}
Let $\mathbf{x}=(x_1,\dots,x_n)$ and $\mathbf{y}=(y_1,\dots,y_n)$ be two sets of variables,
and let $$\Delta(\mathbf{x},\mathbf{y}) = \prod_{i+j \leq n} (x_i - y_j).$$
To each permutation $w\in S_n$ we associate the \emph{double Schubert polynomial}
$$\Sym_w(\mathbf{x},\mathbf{y}) = \partial_{w^{-1} w_0} \Delta(\mathbf{x},\mathbf{y}),$$
where the \emph{divided difference operator} acts on the $x$-variables.
\end{Def}

\begin{Def}
A \emph{partition} $\lambda = (\lambda_1,\dots,\lambda_r)$ is a weakly decreasing sequence of positive integers.
We say that $r$ is the \emph{length} of $\lambda,$
and denote it $r=\length(\lambda).$
\end{Def}

\begin{Def}
The \emph{diagram} or \emph{Rothe diagram} of a permutation $w$ is 
$$D(w) = \{(i,j) \ \vert \ 1 \leq i, j \leq n, w(i)>j, w^{-1}(j)>i\}.$$

The sequence of the numbers of the points of the diagram in successive rows is called
the \emph{Lehmer code} or \emph{code} $c(w)$ of the permutation. We also define $c^{-1}(l)$ to be the permutation whose Lehmer code is $l$. The partition obtained by sorting the components of the code is called the \emph{shape} $\lambda(w)$ of $w$.
\end{Def}

\begin{Ex} If $w = (1,3,5,4,2)$ then $c(w)=(0,1,2,1,0)$ and 
$\lambda(w) = (2,1,1)$.
\end{Ex}

\begin{Def}
We say that a permutation $w$ is \emph{vexillary} if and only if there does not 
exist a sequence $i<j<k<\ell$ such that $w(j)<w(i)<w(\ell)<w(i)$.  Such a permutation
is also called \emph{$2143$-avoiding.}
\end{Def}

\begin{Def}\label{def:flag}
We define the \emph{flag} of a vexillary permutation $w$, starting from its code $c(w)$, in the following fashion.  If $c_i(w) \neq 0$, let $e_i$ be the greatest integer $j\geq i$ such that 
$c_j(w) \geq c_i(w)$.  The flag $\phi(w)$ is then the sequence of integers $e_i$, ordered to be increasing.
\end{Def}

\begin{Def}
Let $X_i$ denote the family of indeterminates $x_1,\dots,x_i$.  
For $d_1,\dots,d_n$ a weakly increasing sequence of $n$ integers,
we define the \emph{flagged Schur function} 
$$s_{\lambda}(X_{d_1},\dots,X_{d_n}) = \sum_T \mathbf{x}^{\type(T)},$$
where the sum runs over the set of semistandard tableaux $T$ with shape $\lambda$ for which 
the entries in the $i$th row are bounded above by $d_i$.
\end{Def}

There is also a notion of flagged double Schur polynomials.  One can define them in terms of 
tableaux or via a Jacobi-Trudi type formula \cite[Section 2.6.5]{Manivel}.

\begin{Th}\cite[Corollary 2.6.10]{Manivel}\label{th:flag}
If $w$ is a vexillary permutation with shape $\lambda(w)$ and with flags $\phi(w) = (f_1,\dots,f_m)$ and $\phi(w^{-1})=(g_1,\dots,g_m)$, then we have
$$\Sym_w(\mathbf{x}; \mathbf{y}) = s_{{\lambda}(w)}(X_{f_1}-Y_{g_m},\dots,X_{f_m}-Y_{g_1}),$$
i.e. the double Schubert polynomial of $w$ is a flagged double Schur polynomial.
\end{Th}

\section{Main results}

Let $w\in S_n$ be a state. In what follows, we write $a\to b \to c$ if the letters $a, b, c$ appear in cyclic order in $w$.
So for example, if $w=1423$, we have that $1\to 2\to 3$ and $2\to 3\to 4$,
but it is not the case that $3 \to 2 \to 1$ or $4\to 3 \to 2$.

\begin{equation}\label{eq:xy}
   xyFact(w)=\prod\limits_{i=1}^{n-2}\prod_{\substack{k>i+1 \\ i \to i+1 \to k}}(x_1-y_{n+1-k})\cdots(x_i-y_{n+1-k}). 
\end{equation}

The following is our main theorem;  when 
each $y_i=0$, it reduces to \cref{thm:main0}.

\begin{Th}\label{thm:main1}
Let $w\in \St(n,k)$, and write  $\Psi(w)=(\lambda^{1},\cdots,\lambda^{k})$. Then the (renormalized) steady state probability is given by  
\begin{equation}\label{eq:maintheorem1}
    \psi_w=xyFact(w)  \prod_{i=1}^k \Sym_{g_n(\lambda^i)},
    \end{equation}
where $\Sym_{g_n(\lambda^i)}$ is the double Schubert polynomial
associated to the permutation with Lehmer code
$g_n(\lambda^i)$, and $g_n$ is given by \cref{def:g}.
\end{Th}

We also prove the \emph{monomial factor conjecture} from 
\cite{LW}.  Suppose that $y_i=0$ for all $i$.
Given a state $w$, let $a_i(w)$ be the number of integers greater than $(i+1)$ on the clockwise path from $(i+1)$ to $i$.  
Let $\eta(w)$ be the largest monomial that can be factored
out of $\psi_w$.    The following statement was conjectured in \cite[Conjecture 2]{LW}.

\begin{Th}
Let $w\in \St(n).$  Then 
$$\eta(w) = \prod_{i=1}^{n-2} x_i^{a_i(w)+\cdots+a_{n-2}(w)}.$$
\end{Th}


\section{Multiline queues and semistandard tableaux}\label{seq:MLQ}

It was proved in \cite{AM} that when each $y_i=0$, the steady state probabilities $\psi_w$ for the 
TASEP on a ring can be expressed in terms of the \emph{multiline queues} of 
 Ferrari and Martin \cite{FM}.  On the other hand, we know from \cref{thm:main0} that 
 when $w\in \St(n,1)$ (i.e. $w^{-1}$ is a Grassmann permutation and $w_1=1$),
 $\psi_w$ equals a monomial times a single flagged Schur polynomial.  In this section we will
 explain that result by giving a bijection between the relevant multiline queues and 
 the corresponding semistandard tableaux.

\begin{Def} Fix positive integers $L$ and $n$. A \emph{multiline queue} $Q$ is an $L \times n$ array in which each of the $Ln$ positions is either vacant or occupied by a ball. 
We say it has \emph{content 
$\mathbf{m} = (m_1,\dots, m_n)$} if it has $m_1+ \dots + m_i$ balls in row 
$i$ for $1 \leq i \leq n$.  We number the rows from top to bottom from 1 to $L$, and the columns from right to left from 1 to $n$. 
\end{Def}

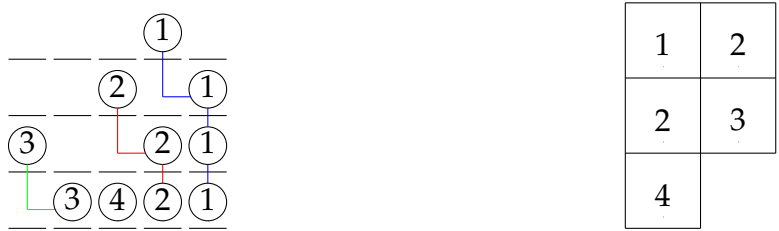
\begin{figure}[h]\centering
\begin{tikzpicture}[scale=0.5]

\draw[-] (-4.6,1) --(-5.6,1);
\draw[-] (-5.8,1) --(-6.8,1);
\draw[-] (-7,1) --(-8,1);
\draw[-] (-8.2,1) --(-9.2,1);
\draw[-] (-9.4,1) --(-10.4,1);

\begin{scope}[shift={(0,1.5)}]

\draw[-] (-4.6,1) --(-5.6,1);
\draw[-] (-5.8,1) --(-6.8,1);
\draw[-] (-7,1) --(-8,1);
\draw[-] (-8.2,1) --(-9.2,1);
\draw[-] (-9.4,1) --(-10.4,1);

\draw (-5.1,1.7) circle (0.5);
\filldraw[black] (-5.1,1.2) circle (0.000001pt) node[anchor=south] {$1$};
\draw (-6.3,1.7) circle (0.5);
\filldraw[black] (-6.3,1.2) circle (0.000001pt) node[anchor=south] {$2$};
\draw (-9.9,1.7) circle (0.5);
\filldraw[black] (-9.9,1.2) circle (0.000001pt) node[anchor=south] {$3$};

\end{scope}

\begin{scope}[shift={(0,3)}]

\draw[-] (-4.6,1) --(-5.6,1);
\draw[-] (-5.8,1) --(-6.8,1);
\draw[-] (-7,1) --(-8,1);
\draw[-] (-8.2,1) --(-9.2,1);
\draw[-] (-9.4,1) --(-10.4,1);

\draw (-5.1,1.7) circle (0.5);
\filldraw[black] (-5.1,1.2) circle (0.000001pt) node[anchor=south] {$1$};

\draw (-7.5,1.7) circle (0.5);
\filldraw[black] (-7.5,1.2) circle (0.000001pt) node[anchor=south] {$2$};

\end{scope}

\begin{scope}[shift={(0,4.5)}]

\draw[-] (-4.6,1) --(-5.6,1);
\draw[-] (-5.8,1) --(-6.8,1);
\draw[-] (-7,1) --(-8,1);
\draw[-] (-8.2,1) --(-9.2,1);
\draw[-] (-9.4,1) --(-10.4,1);

\draw (-6.3,1.7) circle (0.5);
\filldraw[black] (-6.3,1.2) circle (0.000001pt) node[anchor=south] {$1$};
\end{scope}

\draw[red] (-7.5,4.2)-- (-7.5,3);
\draw[red] (-7.5,3)-- (-6.75,3);
\draw[red] (-6.3,2.7)-- (-6.3,2.2);

\draw[blue] (-7.5+1.2,4.2+1.5)-- (-7.5+1.2,3+1.5);
\draw[blue] (-7.5+1.2,3+1.5)-- (-6.75+1.2,3+1.5);
\draw[blue] (-5.1,4.2)-- (-5.1,3.7);
\draw[blue] (-5.1,2.7)-- (-5.1,2.2);

\draw[green] (-7.5-2.4,4.2-1.5)-- (-7.5-2.4,3-1.5);
\draw[green] (-7.5-2.4,3-1.5)-- (-6.75-2.4,3-1.5);

\draw (-5.1,1.7) circle (0.5);
\filldraw[black] (-5.1,1.2) circle (0.000001pt) node[anchor=south] {$1$};
\draw (-6.3,1.7) circle (0.5);
\filldraw[black] (-6.3,1.2) circle (0.000001pt) node[anchor=south] {$2$};
\draw (-8.7,1.7) circle (0.5);
\filldraw[black] (-8.7,1.2) circle (0.000001pt) node[anchor=south] {$3$};
\draw (-7.5,1.7) circle (0.5);
\filldraw[black] (-7.5,1.2) circle (0.000001pt) node[anchor=south] {$4$};

\draw[-] (6,1) --(6,7);
\draw[-] (6,7) --(10,7);
\draw[-] (6,5) --(10,5);
\draw[-] (10,7) --(10,3);
\draw[-] (10,3) --(6,3);
\draw[-] (8,7) --(8,1);
\draw[-] (6,1) --(8,1);

\filldraw[black] (7,1.3) circle (0.000001pt) node[anchor=south] {$4$};
\filldraw[black] (7,3.3) circle (0.000001pt) node[anchor=south] {$2$};
\filldraw[black] (7,5.3) circle (0.000001pt) node[anchor=south] {$1$};
\filldraw[black] (9,3.3) circle (0.000001pt) node[anchor=south] {$3$};
\filldraw[black] (9,5.3) circle (0.000001pt) node[anchor=south] {$2$};

\end{tikzpicture}
\caption{A multiline queue of type $(1,2,4,3,5)$, and the corresponding semistandard tableau under the bijection in \cref{Prop:bijection}.} \label{fig0}
\end{figure}

\begin{Def}
Given an $L\times n$ multiline queue $Q$, the \emph{bully path projection} on $Q$ is, for each row $r$ with $1 \leq r \leq L-1$, a particular matching of balls from row $r$ to row $r+1$, which we now define. If ball $b$ is matched to  ball $b_0$ in the row below then we connect $b$ and $b_0$ by the shortest path that travels either straight down or from left to right (allowing the path to wrap around the cylinder if necessary). 
Here each ball is assigned a \emph{class}, and matched according to the following algorithm:
\begin{itemize}
    \item All the balls in the first row are defined to be of class $1$.
    \item Suppose we have matched all the balls in rows $1,2,\dots, r-1$ and have assigned a class to all balls in rows $1,2,\dots,r$.  We now consider the balls in rows $r$.
    \item Pick any order of the balls in row $r$ such that balls with smaller labels come before balls with larger labels.  Consider the balls in this order;   
    suppose we are considering a ball $b$ of class $i$ in row $r$. If there is an 
    unmatched ball directly below
    $b$ in row $r+1$, we let $M(b)$ be that ball; otherwise we move to the right in row $r+1$ and let $M(b)$ be the first unmatched ball that we find (wrapping around from column $1$ to $n$ if necessary).  We match $b$ to ball $M(b)$ and say that $M(b)$ is of class $i$.  
    \item The previous step gives a matching 
    of all balls in row $r$ to balls below in row $r+1$.   We assign class $r+1$ to any 
    balls in row $r+1$ that were not yet assigned a class.  We now repeat the process 
    and consider the balls in row $r+1$.
    \end{itemize}

After completing the bully path projection for $Q$, let $w=(w_1,\cdots,w_n)$ be the labeling of the balls read from right to the left in row $L$ (where a vacancy is denoted by $L+1$). We say that $Q$ is a multiline queue of \emph{type $w$} and let $MLQ(w)$ denote the set of all multiline queues of type $w$.  We also consider a type of row $r$ in $Q$ to be the labeling of the balls read from right to the left in row $r$ (where a vacancy is denoted by $r+1$).

 A vacancy in $Q$ is called $i-covered$ if it is traversed by a path starting on row $i$, but not traversed by any path starting on row $i'$ such that $i'<i$.
\end{Def}
See \cref{fig0} for an example.

We define a weight $wt(Q)$ for multiline queues. It was first introduced in \cite{AL}.
\begin{Def}
Given an $L\times n$ multiline queue $Q$, let $v_r$ be the number of vacancies in row $r$ and let $z_{r,i}$ be the number of $i-covered$ vacancies in row $r$. 
Set $V_i=\sum\limits_{j=i+1}^{L}v_j$. We define 
\begin{equation*}
    wt(Q)=\prod\limits_{i=1}^{L-1}(x^{V_i}_i)\prod\limits_{1\leq i<r\leq L} (\frac{x_r}{x_i})^{z_{r,i}}.
\end{equation*}
\end{Def}
\begin{Ex}
The multiline queue $Q$ in \cref{fig0} has a $1-covered$ vacancy in row 2, a $2-covered$ vacancy in row 3 and a $3-covered$ vacancy in row 4. The weight of $Q$ is 
\begin{equation*}
    wt(Q)=x^{3+2+1}_1 x^{2+1}_2 x^{1}_3 (\frac{x_2}{x_1}) (\frac{x_3}{x_2})(\frac{x_4}{x_3})=x^{5}_1x^{3}_2x_3x_4.
\end{equation*}
\end{Ex}

The following result was conjectured in \cite{AL} and proved in \cite{AM}.
\begin{Th}\cite{AM}
Consider the inhomogeneous TASEP on a ring (with each $y_i=0$).  We have 
\begin{equation*}
\psi_w=\sum\limits_{Q\in MLQ(w)}wt(Q).    
\end{equation*}
\end{Th}

We now give a (weight-preserving up to a constant factor) bijection between multiline queues in $MLQ(w)$ and certain semistandard tableaux, when $w\in \St(n,1)$, i.e. $w^{-1}$ is a Grassmann permutation and $w_1=1$.

\begin{figure}[h]\centering
\begin{tikzpicture}[scale=0.6]

\draw[-] (0,0) --(1,0);
\draw[-] (0,1) --(2,1);
\draw[-] (0,2) --(2,2);
\draw[-] (0,3) --(2,3);

\draw[-] (0,0) --(0,3);
\draw[-] (1,0) --(1,3);

\draw[-,blue,thick] (2,3) --(2,1);
\draw[-,blue,thick] (2,1) --(1,1);
\draw[-,blue,thick] (1,1) --(1,0);
\draw[-,blue,thick] (1,0) --(0,0);

\filldraw[black] (2,2.5) circle (0.000001pt) node[anchor=west] {1};
\filldraw[black] (2,1.5) circle (0.000001pt) node[anchor=west] {2};
\filldraw[black] (1,0.5) circle (0.000001pt) node[anchor=west] {3};

\filldraw[red] (1.5,1) circle (0.000001pt) node[anchor=south] {4};
\filldraw[red] (0.5,0) circle (0.000001pt) node[anchor=south] {5};

\end{tikzpicture}
\caption{The partition $\lambda=(2,2,1)$ and $w(\lambda)=(1,2,4,3,5).$} \label{fig5}
\end{figure}
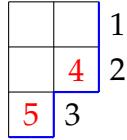

\begin{Def} \label{def:wlambda} Given a partition $\lambda=(\mu_1^{b_1},\cdots,\mu_k^{b_k},0^{c})$, such that $\mu_1>\cdots>\mu_k>0$ and $b_i>0, c\geq0$, we 
define a permutation $w(\lambda)$ as follows. Identify $\lambda$ with the lattice path from $(\mu_1,\sum\limits_{i=1}^{k}b_i+c)$ to $(0,0)$ that defines the southeast border of its Young diagram. 
Label the vertical steps of the lattice path from $1$ to $k$ from top to bottom, and then the horizontal steps in increasing order from right to left starting from $k+1$.  Reading off the numbers along the lattice path gives $w(\lambda)$. See \cref{fig5}.
\end{Def}



\begin{Prop}\label{Prop:bijection}
Given a partition $\lambda=(\mu_1^{b_1},\cdots,\mu_k^{b_k},0^{c})$ as in \cref{def:wlambda}, 
let $d=(d_1,\cdots,d_k)$ be the numbers assigned to horizontal steps right after vertical steps in the construction of $w(\lambda)$. For example, in \cref{fig5}, $d=(4,5)$. Let $d'$ be the vector 
\begin{equation*}
    d'=(\underbrace{d_1-b_1,\cdots,d_1-1}_\text{$b_1$}, \underbrace{d_2-b_2,\cdots,d_2-1}_\text{$b_2$},\dots,
\underbrace{d_k-b_k,\cdots,d_k-1}_\text{$b_k$}).
\end{equation*}

Then there exists a bijection $f: MLQ(w) \rightarrow SSYT(\lambda,d')$ such that 
$wt(Q)=K x^{type(f(Q))}$ for some monomial $K$, where $SSYT(\lambda,d')$ is the set of semistandard tableaux with shape $\lambda$ for which the entries in the $i$ th row are bounded above by $ d'_i$. In particular, we have 
\begin{equation*}
    \psi_{w(\lambda)}=\sum\limits_{Q\in MLQ(w(\lambda))}wt(Q)=K\sum\limits_{T\in SSYT(\lambda,d')}x^{type(T)}=K s_{\lambda}(X_{d'_1}, X_{d'_2},\dots).
\end{equation*}

\end{Prop}

\printbibliography

\end{document}